\documentclass[12pt]{amsart}
\usepackage{graphicx}
\usepackage{amsmath}

\newtheorem{theorem}{Theorem}[section]
\theoremstyle{plain}

\newtheorem{corollary}{Corollary}[section]

\newtheorem{lemma}{Lemma}[section]

\numberwithin{equation}{section}

\input{epsf}

\begin{document}
\title[Markov Processes on Graphs]
{Continuous Time Markov Processes on Graphs}
\author{Jianjun Tian}
\address[]{Mathematical Biosciences Institute\\ 
The Ohio State University\\
Columbus, OH 43210, USA}
\email[]{tianjj@mbi.ohio-state.edu}
\author{Xiao-Song Lin}
\address[]{Department of Mathematics\\
University of California\\ Riverside, CA 92521, USA}
\email[]{xl@math.ucr.edu}
\begin{abstract} We study continuous time Markov processes on graphs. The notion of frequency is introduced, which serves well as a scaling factor between any Markov time of a continuous time Markov process and that of its jump chain. As an application, we study ``multi-person simple random walks'' on a graph $G$ with $n$ vertices. There are $n$ persons distributed randomly at the vertices of $G$. In each step of this discrete time Markov process, we randomly pick up a person and move it to a random adjacent vertex. We give estimate on the expected number of steps for these $n$ persons to meet all together at a specific vertex, given that they are at different vertices at the begininng. For regular graphs, our estimate is exact.   
\end{abstract}
\maketitle

\section{Introduction}

In this paper, for simplicity, we will consider connected simple graphs only. These are connected graphs without multiple edges and loops. We will adopt the following notations and terminologies for a graph $G$. The sets of vertices and edges of $G$ are 
$V(G)$ and $E(G)$, respectively. The {\it order} $n$ of $G$ is the number of vertices of $G$, and the {\it size} $m$ of $G$ is the number of edges of $G$. 
Thus, $n=|V(G)|$ and $m=|E(G)|$. For a vertex $x\in V(G)$, $\Gamma(x)$ is the set of vertices which are connected to $x$ by an edge in $E(G)$. The {\it degree} of a vertex $x$ is $d(x)=|\Gamma(x)|$. We have
$$\sum_{x\in V(G)}\,d(x)=2m.$$
The adjacent matrix of $G$ is denote by $A(G)$ and the diagonal matrix $D(G)$ has the sequence of degrees at each vertex as its diagonal entries.
Finally, we denote
$$d_m=\min\,\{d(x)\,;\,x\in V(G)\}\quad\text{and}\quad
d_M=\max\,\{d(x)\,;\,x\in V(G)\}.$$
\smallskip

What concerns us primarily in this paper is the following continuous time Markov process on a graph $G$: The probability that a person standing at a vertex $x$ of $G$ to jump to a neighboring vertex in $\Gamma(x)$ within a time period $\Delta t$ is $d(x)\Delta t+o(\Delta t)$, and once jumped, the person at $x$ has equal probability to land at a vertex $y\in\Gamma(x)$. If we write 
$$Q=Q(G)=-D(G)+A(G),$$
the transition probability matrix of this Markov process is
$$P(t)=e^{tQ}.$$
We call this Markov process CTSRW (continuous time simple random walks). In the literature, it is the discrete time simple random walks (SRW) on a graph $G$ that people concern most. One usually call SRW the {\it jump chain} of CTSRW. The transition probability matrix of SRW is $D(G)^{-1}A(G)$.  
\smallskip

We introduce in this paper a fundamental quantity for CTSRW on a graph $G$ called {\it frequency}. Let $N(t)$ be the expected number of jumps of the Markov process CTSRW up to time $t$. Then the frequency $f$ of CTSRW is defined to be 
$$f=\lim_{t\rightarrow\infty}\,\frac{N(t)}{t}.$$
Using the L\'evy formula, we are able to calculate the frequency for CTSRW and get
$$f=\frac{2m}{n}.$$
\smallskip

The frequency turns out to be a natural scaling factor between various important quantities of CTSRW and SRW, respectively. For example, we have the following theorem.

\begin{theorem} Let $G$ be non-bipartite. For a vertex $x$ of $G$, let $T_x$ be the first return time of CTSRW on $G$ and $N_{T_x}$ be the number of jumps during the time period $[0,T_x]$. Notice that $N_{T_x}$ is the first return time for the discrete time SRW. Then, the expectations of $T_x$ and $N_{T_x}$, $E(T_x)$ and $E(N_{T_x})$, respectively, are related by the following equation:
$$E(N_{T_x})=fE(T_x).$$    
\end{theorem}
 
More generally, we have the following theorem.

\begin{theorem} Let $T$ be any Markov time (or stopping time) of CTSRW on a graph $G$ with finite expectation, let $N_T$ be the number of jumps during the time period $[0,T]$. Then
$$d_m\leq f\leq d_M$$
and
$$d_mE(T)\leq E(N_T)\leq d_M E(T).$$
In particular, if $G$ is regular so that $d_m=d_M$, we have $E(N_T)=fE(T)$.
\end{theorem}

As an application, we consider {\it multi-person simple random walks} (MPSRW) on a graph $G$. To start with, we assume that each of the $n$ vertices of $G$ is occupied by a person. In each step of this Markov chain, there is one person, equally possible for each of these $n$ persons, who moves to a neighboring vertex, also equally possible for each of the neighboring vertices. We would like to know the expected number of steps this Markov chain should take for these $n$ persons to meet all together at a specified vertex. 
\smallskip

We will see that this Markov chain is the jump chain of a continuous time Markov process on the set $M_n$ of {\it maps} from $\{1,2,\dots,n\}$ to itself.
On the other hand, this continuous time Markov process on $M_n$ can be identified with the $n$-th tensor power of CTSRW on $G$. Thus, computation of expectations of various stopping times for this continuous time Markov process on $M_n$ can be carried out. We are then able to use Theorem 1.2 mentioned above to obtain estimates for the expected number of steps for MPSRW on $G$ to coalesce.  
\smallskip

We refer the reader to \cite{Bo, GR, Lo} for basic terminologies and results in
the study of simple random walks on graphs.

\section{Continuous time Markov process on weighted graphs}

Let $G$ be a connected weighted graph with order $n$ size $m$, we define
continuous time Markov process $X_t$ on $G$ by giving its infinitesimal generator 
$Q$ as the negative weighted Laplacian of $G$. Specifically, every edge 
$xy\in E(G)$ is associated with a positive number (weight) $w_{xy}$. We do not direct edges in $G$ and therefore $w_{xy}=w_{yx}$. We denote by
$$w_x=\sum_{y\in\Gamma(x)}\,w_{xy}$$
the total weight at the vertex $x$. We write $Q=\left(q_{xy}\right) _{n\times n}$, where  
$$
q_{xy}=\begin{cases}
w_{xy}&\text{if $xy\in E(G)$,} \\ 
-w_{x}&\text{if $x=y$,} \\ 
0&\text{otherwise.}
\end{cases}
$$
Thus the probability transition matrix $X_t$ is
given by 
\begin{equation*}
P\left( t\right) =e^{tQ}=\left( p_{xy}\left( t\right) \right) _{n\times n},
\end{equation*}
and transition probability from vertex $x$ to vertex $y$ is given by 
\begin{equation*}
\Pr \left\{ X\left( h+t\right) =y\mid X\left( h\right) =x\right\}
=p_{xy}\left( t\right) .
\end{equation*}
By the way, we may call $-Q=L_{w}$ the weighted Laplacian of the weighted
graph $G.$
\smallskip

In the special case of $w_{xy}=1$ for all $xy\in E(G)$, we have CTSRW on the graph $G$. The infinitesimal generator 
$Q=-L(G)=-D(G) +A(G)=(q_{xy})_{n\times n}$ is given by 
$$
q_{xy}=\begin{cases}
1 &\text{if $xy\in E(G)$,} \\ 
-d(x) &\text{if $x=y$,} \\ 
0& \text{otherwise.}
\end{cases}
$$

\subsection{Ergodicity}
We set $u=\left( \frac{1}{n},\frac{1}{n},\cdots ,\frac{1}{n}\right) $ to be
a probability vector. Then, since $Q$ is symmetric, $u$ is an invariant measure of the Markov process $X_t$. Namely, we have $uQ=0$ and 
\begin{equation*}
uP\left( t\right) =u\left( I+tQ+\frac{t^{2}}{2!}Q^{2}+\cdots \right) =u.
\end{equation*}
We claim that $u$ also is the ergodic vector or the stationary distribution.
To see this, notice first that the graph $G$ is connected so the process $X_t$ is irreducible. Thus 
$$\lim_{t\rightarrow \infty }p_{xy}\left( t\right) =v_{y}$$
exists and does not depend on $x.$ 
Actually, we have $v_{y}=u_{y}=\frac1{n}$ by the following calculation:
First, for any fixed $t>0$,
$$\begin{aligned}
uP(2t) &=uP\left( t\right) P\left( t\right) =uP\left( t\right) =u, \\
&\dots \dots \\
uP\left( kt\right) &=uP\left( \left( k-1\right) t\right) =\cdots =uP(t)=u.
\end{aligned}
$$
Then
\begin{equation*}
u_{y}=\sum_{x}u_{x}p_{xy}\left( kt\right) ,\text{ \ for \ }y\in V\left(
G\right).
\end{equation*}
Let $k\rightarrow \infty$, we get 
\begin{equation*}
u_{y}=\sum_{x}u_{x}v_{y}=v_{y}.
\end{equation*}

\subsection{The mean first return time}
For a vertex $x\in V(G)$, denote by $T_{xx}$ the first return time to $x$, given that the Markov process $X_t$ starts at $x$. That is 
$$T_{xx}=\inf\,\{ t\,:\,t>\rho _{x},X_t=x\mid X_0=x\}$$ 
where $\rho_{x}$ is the exit time from the vertex $x$. We denote by $h(x,x)$
the mean first return time $E(T_{xx})$.

\begin{lemma} The
mean first return time is $$h(x,x)=\frac{n}{w_{x}}.$$ 
\end{lemma}

\begin{proof}
Write 
\begin{equation*}
F_{xx}(t)=\Pr\,\{ T_{xx}\leq t\}. 
\end{equation*}
Then we have the equation \cite{Ch}
\begin{equation*}
p_{xx}=e^{-w_xt}+\int_{0}^{t}p_{xx}\left( t-s\right)
dF_{xx}\left( s\right). 
\end{equation*}
Taking Laplace transform, we get
$$\begin{aligned}
\phi_{xx}\left( \lambda \right)  &=\int_{0}^{\infty }e^{-\lambda
t}p_{xx}\left( t\right) dt\\
&=\int_{0}^{\infty }e^{-\left( \lambda +w_x\right) t}dt+\int_{0}^{\infty }e^{-\lambda
t}\int_{0}^{t}p_{xx}\left( t-s\right) dF_{xx}\left( s\right) dt \\
&=\frac{1}{\lambda +w_x}+\int_{0}^{\infty }\int_{s}^{\infty
}e^{-\lambda t}p_{xx}\left( t-s\right) dtdF_{xx}\left( s\right)  \\
&=\frac{1}{\lambda +w_x }+\int_{0}^{\infty }\int_{0}^{\infty
}e^{-\lambda \left( v+s\right) }p_{xx}\left( v\right) dvdF_{xx}\left(
s\right)  \\
&=\frac{1}{\lambda +w_x}+\int_{0}^{\infty }e^{-\lambda
s}\int_{0}^{\infty }e^{-\lambda v}p_{xx}\left( v\right) dvdF_{xx}\left(
s\right)  \\
&=\frac{1}{\lambda +w_x}+\int_{0}^{\infty }e^{-\lambda s}\phi
_{xx}\left( \lambda \right) dF_{xx}\left( s\right)  \\
&=\frac{1}{\lambda +w_x}+\phi _{xx}\left( \lambda \right)
l_{x}\left( \lambda \right) 
\end{aligned}
$$
where $l_{x}\left( \lambda \right) =\int_{0}^{\infty }e^{-\lambda
s}dF_{xx}(s).$ Then, we have 
\begin{equation*}
\lambda \phi _{xx}\left( \lambda \right) =\frac{1}{\lambda +d\left( x\right) 
}\left( \frac{1-l_{x}\left( \lambda \right) }{\lambda }\right) ^{-1}.
\end{equation*}
\smallskip
Since
\begin{equation*}
\lim_{\lambda \rightarrow 0}\frac{1-l_{x}\left( \lambda \right) }{\lambda }%
=\lim_{\lambda \rightarrow 0}\int_{0}^{\infty }se^{-\lambda s}dF_{xx}\left(
s\right) =\lim_{\lambda \rightarrow 0}\int_{0}^{\infty }sdF_{xx}\left(
s\right) =h\left( x,x\right) 
\end{equation*}
and 
\begin{equation*}
\lim_{\lambda \rightarrow 0}\lambda \phi_{xx}\left( \lambda \right) =u_{x},
\end{equation*}
we get 
\begin{eqnarray*}
u_{x} &=&\lim_{\lambda \rightarrow 0}\frac{1}{\lambda +w_x}
\left( \frac{1-l_{x}\left( \lambda \right) }{\lambda }\right) ^{-1} \\
&=&\frac{1}{h\left( x,x\right) d\left( x\right) }.
\end{eqnarray*}
Thus,
$$h(x,x)=\frac{n}{w_x}.$$
\end{proof}

\subsection{The mean hitting time of $y$ from $x$}
In general, define
\begin{equation*}
T_{xy}=\inf\,\{ t\,:\,t>\rho _{x},X_t=y\mid X_0=x\}.
\end{equation*}
I. e. $T_{xy}$ is the time of first entrance into,
or hitting, the vertex $y$, given that the process starts at $x$.
We denote the mean first hitting time of $y$ from $x$
by $h\left(x,y\right)=E(T_{xy})$. We have the following equation:
\begin{eqnarray*}
h(x,x)=\frac{1}{w_{x}}+\sum_{y\in \Gamma \left( x\right) }\frac{w_{xy}}{%
w_{x}}\,h\left( y,x\right) 
=\frac{1}{w_{x}}+\frac{1}{w_{x}}\sum_{y\in \Gamma \left( x\right)
}w_{xy}\,h(y,x).
\end{eqnarray*}

\begin{lemma}
The mean hitting time can be calculated from the following integral:  
\begin{equation*}
h(x,y)=n\int_{0}^{\infty }\left( p_{yy}\left( t\right) -p_{xy}\left(
t\right) \right) dt.
\end{equation*}
\end{lemma}

The formula in this lemma is similar to the formula in the
discrete time case. We omit the proof since it is also analogous to the 
discrete time case. 

\subsection{The stationary distribution of SRW}
For an unweighted graph $G,$ the jump chain CTSRW on $G$ is SRW on $G$. 
We know that $Q=-D+A,$ and the transition probability matrix of SRW
is $D^{-1}A$. If we set $\widetilde{\pi}=uD$, where $u$ is the
stationary distribution of CTSRW, then $
uQ=-uD+uA=0$.
So, substitute $u=\ \widetilde{\pi }D^{-1}$ we have 
\begin{equation*}
\ \widetilde{\pi }=\ \widetilde{\pi }D^{-1}A.
\end{equation*}
Thus, $\widetilde{\pi}$ is an invariant measure of SRW. We need to
normalize it. Let
$$\pi =\frac{1}{\sum_{i=1}^{n}\ \widetilde{\pi }_{i}}\ 
\widetilde{\pi }.$$ 
Then $\pi$ is an invariant distribution for SRW. Specifically 
\begin{equation*}
\pi _{x}=\frac{d\left( x\right) }{\sum_{x\in V(G)}d(x)}=\frac{d(x)%
}{2m}.
\end{equation*}
\smallskip

If the graph $G$ is non-bipartite, this invariant distribution is also
the stationary distribution. It is well known that the mean number of steps SRW should take to return to the vertex $x$ for the first time is 
$1/\pi_x=2m/d(x)$. 
Recall that the mean first return time of CTSRW is $n/d(x)$. Therefore it is natural to think of the quantity
$$f=\frac{2m/d(x)}{n/d(x)}=\frac{2m}{n}$$
as the {\it frequency} (number of jumps per unit time) of CTSRW.  We will make this notion precise in the following subsection. 

\subsection{The frequency} To define the frequency for the continuous time Markov process $X_t$ on a weighted graph, we first define a quantity $N(t)$ for $t>0$:
\begin{equation*}
N\left( t\right) =E\left( \text{the number of jumps of $X_{t}$ up to
time }t\right).
\end{equation*}

\begin{theorem} We have  
\begin{equation*}
f:=\lim_{t\rightarrow \infty }\frac{N\left( t\right)}{t}=\frac{2w}{n}
\end{equation*}
where $w$ is the total weight of $G$
$$w=\frac12\,\sum_{xy\in E(G)}\,w_{xy}.$$
\end{theorem}

\begin{proof} Let us recall the L\'{e}vy formula first. See \cite{Ch}. Given a Markov process $X_{t}$, we consider a purely
discontinuous functional \ $A=\left\{ A_{t}:0<t< \infty \right\} $ on the 
path space defined by 
\begin{equation*}
A_{t}=\sum_{0< s\leq t}g\left( X_{s^{-}},X_{s}\right),\qquad t> 0,
\end{equation*}
where $g$ is a function on $V(G)\times V(G)$.
Also, we define a function $b_{Q}$ on $V(G)$ by 
\begin{equation*}
b_{Q}\left( x\right) =\sum_{y\neq x}w_{xy}g\left( x,y\right) ,\text{ \ }x\in
V(G),
\end{equation*}
and the integral functional $B=\left\{ B_{t}:0\leq t< \infty \right\} 
$ on the path space is defined by  
\begin{equation*}
B_{t}=\int_{0}^{t}b_{Q}\left( X_{s}\right) ds=\int_{0}^{t}\sum_{y\neq
X_{s}}w_{X_{s},y}g\left( X_{s},y\right)ds .
\end{equation*}
Then, the relationship between the functionals $A$ and $B$ is given by the
L\'{e}vy formula: 
\begin{equation*}
E_{x}\sum_{0< s\leq t}g\left( X_{s^{-}},X_{s}\right) \alpha \left(
t\right) =E_{x}\int_{0}^{t}\alpha \left( t\right) b_{Q}\left( X_{s}\right)
ds,\qquad t>0,
\end{equation*}
for any continuous positive function $\alpha \left( t\right) .$

Now, taking $\alpha \left( t\right) =1,$ the L\'{e}vy formula tells us 
\begin{equation*}
E_{x}A_{t}=E_{x}B_{t}=\int_{0}^{t}P\left( s\right) \cdot b_{Q}\left(
x\right) ds.
\end{equation*}
Furthermore, let $$g\left( x,y\right) =\left\{ 
\begin{array}{c}
1,\text{ if \ }x\neq y, \\ 
0,\text{ if \ }x=y.
\end{array}
\right. $$
Then $A_{t}$ is the number of transitions of states of $X_{t}$ up to
time $t$, i. e. $EA_{t}=N\left( t\right)$.  
\smallskip

We start at the vertex $x$. Then 
\begin{equation*}
E_{x}A_{t}=\int_{0}^{t}P\left( s\right) \cdot b_{Q}\left( x\right)
ds=\int_{0}^{t}\sum_{y\in V(G)}p_{xy}\left( s\right) w_{y}ds.
\end{equation*}
If we start at an initial distribution $\theta $ on graph $G,$ then 
\begin{equation*}
E_{\theta }A_{t}=\int_{0}^{t}\theta P\left( s\right) \cdot
b_{Q}ds=\int_{0}^{t}\sum_{x,y}\theta _{x}p_{xy}\left( s\right) w_{y}ds.
\end{equation*}
Thus, we have 
$$\begin{aligned}
\lim_{t\rightarrow \infty }\frac{E_{\theta }A_{t}}{t}&=\lim_{t\rightarrow
\infty }\frac{\int_{0}^{t}\sum_{x,y}\theta _{x}p_{xy}\left( s\right) w_{y}ds
}{t}\\
&=\lim_{t\rightarrow \infty }\sum_{x,y}\theta _{x}p_{xy}\left( t\right)
w_{y}\\
&=\sum_{x,y}\theta _{x}u_{y}w_{y}=\sum_{y}u_{y}w_{y}
\end{aligned}
$$
and it is independent of the initial condition $\theta$. So, 
\begin{equation*}
f=\sum_{y\in V(G)}\frac{1}{n}\,w_{y}=\frac{2w}{n}.
\end{equation*}
\end{proof}

Using the notion of frequency, we can compare various Markov times for the continuous time Markov process and its jump chain. Let us recall the notion of Markov time (or stopping time) first. Associated with a stochastic process, there are random variables independent of the future. This kind of random variables are called
Markov time or stopping time. Specifically, let $\sigma$ be a non-negative
random variable associated with a given process 
$\left\{ X_{t}:0\leq t\leq \infty \right\} $. In the other words, $\sigma$
associates with each sample function $X_{t}$ a nonnegative number which we
denote by $\sigma \left( X_{t}\right)$. Such a random variable $\sigma $ is said to be a Markov
time relative to the process $X_{t}$ if it has the following property:
\smallskip

If $X_{t}$ and $Y_{t}$ are two sample functions of the process such that
$X_\tau =Y_{\tau }$ for $0\leq \tau \leq s$ and $\sigma \left(
X_{t}\right)< s,$ then $\sigma \left( X_{t}\right) =\sigma \left(
Y_{t}\right)$.
\smallskip

Now, let's state our main result.

\begin{theorem}
Let $T$ be any Markov time with finite expectation, i.e. $E(T)<
\infty$ associated with the continuous Markov process $X_t$ on a weighted graph $G$. Let $N_{T}$ be the number of jumps of $X_{t}$ during the period $[0,T]$.
Then, we have 
$$w_{m} \leq f\leq w_{M}$$
and
$$w_{m}E(T) \leq E(N_{T})\leq w_{M}E\left( T\right),$$
where $w_{m}=\min\,\{w_x\,;\,x\in V(G)\}$ and $w_{M}=\max\,\{w_x\,;\,x\in V(G)\}$.
\end{theorem}

\begin{proof}
By the L\'{e}vy formula, for an initial distribution $\theta$, we
have 
\begin{equation*}
N\left( t\right) =E_{\theta }A_{t}=\int_{0}^{t}\sum_{x,y\in V(G)}\theta
_{x}p_{xy}\left( s\right) w_{y}ds.
\end{equation*}
Since
\begin{equation*}
\sum_{x,y\in V(G)}\theta _{x}p_{xy}\left( s\right) w_{y}\leq \sum_{x,y\in
V(G)}\theta _{x}p_{xy}\left( s\right) w_{M}=w_{M}
\end{equation*}
and 
\begin{equation*}
\sum_{x,y\in V(G)}\theta _{x}p_{xy}\left( s\right) w_{y}\geq \sum_{x,y\in
V(G)}\theta _{x}p_{xy}\left( s\right) w_{m}=w_{m},
\end{equation*}
we get 
\begin{equation*}
w_{m}t\leq N(t)\leq w_{M}t.
\end{equation*}
This is 
\begin{equation*}
w_{m}\leq \frac{N(t)}{t}\leq w_{M}.
\end{equation*}
By taking limit, we have 
\begin{equation*}
w_{m}\leq f\leq w_{M}.
\end{equation*}
\smallskip

We suppose that Markov time $T$ has distribution \ $F(t)=\Pr \left\{ T<
t\right\}$. Then we have
\begin{equation*}
\int_{0}^{\infty }w_{m}tdF\leq \int_{0}^{\infty }\frac{N(t)}{t}tdF\leq
\int_{0}^{\infty }w_{M}tdF.
\end{equation*}
This actually is 
\begin{equation*}
w_{m}\int_{0}^{\infty }tdF\leq \int_{0}^{\infty }N(t)dF\leq
w_{M}\int_{0}^{\infty }tdF.
\end{equation*}
We recall the conditional expectation 
\begin{equation*}
E(N_{T})=E\left( E\left( A_{T}\mid T=t\right) \right) =\int_{0}^{\infty
}N\left( t\right) dF.
\end{equation*}
Therefore, we get 
\begin{equation*}
w_{m}E(T)\leq E(N_{T})\leq w_{M}E\left( T\right) .
\end{equation*}
\end{proof}

The following are two interesting corollaries. The proofs of them are obvious, so we just state the results.

\begin{corollary} For CTSRW on a graph $G$ and any Markov time $T$ with finite expectation, we have
$$d_{m}\leq f\leq d_{M}$$
and
$$d_{m}E(T) \leq E(N_{T})\leq d_{M}E(T).$$
\end{corollary}

\begin{corollary}
If $G$ is a regular graph with constant degree $d$ at each vertex, then $f=d$  and $E(N_{T})=fE(T)$ for any Markov time $T$ associated with CTSRW on $G$.
\end{corollary}

We may call the inequality in Theorem 2.3 and Corollary 2.1 ``time-step inequality''. Of
course, we have another version as 
\begin{equation*}
\frac{E(N_{T})}{d_{M}}\leq E(T)\leq \frac{E(N_{T})}{d_{m}}.
\end{equation*}
In a sense, those inequalities characterize the timing difference between CTSRW and SRW on a graph. 
It is also interesting to see that the frequency of an unweighted graph is
the average of the eigenvalues of its Laplacian. Let $\lambda _{1}<
\lambda _{2}\leq \lambda _{3}\leq \cdots \leq \lambda _{n}$ be the spectrum
of the Laplacian $L(G)$ of a graph $G$, $L\left( G\right) =D\left( G\right) -A\left(
G\right)$. Then, $\sum_{i=1}^{n}\lambda _{i}=\sum_{x\in V(G)}d\left(
x\right) =2m.$ Thus 
\begin{equation*}
f=\frac{\sum_{i=1}^{n}\lambda _{i}}{n}.
\end{equation*}

\begin{theorem}
Let $G$ be a non-bipartite graph, and $T_{x}$ be the
first return time of CTSRW on $G$. Then $E(N_{T_{x}})=fE(T_{x}).$
\end{theorem}

\begin{proof}
By the theory of discrete time simple random
walks on a graph $G,$ we know $E(N_{T_{x}})=\frac{2m}{d\left( x\right) }$.
\smallskip

For CTSRW on $G$, we know $E(T_{x})=\frac{n}{d\left( x\right) },$ and 
also $f=\frac{2m}{n}$. Thus 
\begin{equation*}
E(N_{T_{x}})=\frac{2m}{d\left( x\right) }=\frac{n}{d\left( x\right) }\cdot 
\frac{2m}{n}=fE(T_{x}).
\end{equation*}
\end{proof}

Now, we consider a special problem as that in SRW. Let $G$ be a connected
non-bipartite graph. We start our CTSRW at a vertex $x$ and fix a neighboring vertex $y$ of $x$.
What is the expected time that our
CTSRW should take in order to return to $x$ through the edge $yx$?
\smallskip

For SRW on $G$, we know the corresponding quantity, the expected number of steps
one should take in order to return to $x$ through
the edge $yx$, is $2m$ [Bollob\'{a}s'].
\smallskip

To deal with the problem for CTSRW on $G$, we formulate the following Markov time:
\begin{equation*}
T_{\left( x,\overline{yx}\right) }=\inf \left\{ t+\rho _{y}:t> 0,X\left(
s+t\right) =y,X\left( s+t+\rho _{y}\right) =x\mid X\left( s\right)
=x\right\}.
\end{equation*}
This is the fixed edge first return time, then $E(N_{T_{\left( x,\overline{yx}\right)
}})=2m.$ By our time-step inequality, we have the following corollary.

\begin{corollary}
The mean fixed edge first return time has bounds as 
\begin{equation*}
\frac{2m}{d_{M}}\leq E\left( T_{\left( x,\overline{yx}\right) }\right) \leq 
\frac{2m}{d_{m}}.
\end{equation*}
\end{corollary}

The following is another case where the frequency gives us a perfect scaling factor between corresponding quantities of CTSRW and SRW, respectively.

\begin{lemma}
Let $G$ be a connected graph of order $n$ and size $m.$ The mean hitting
time $h\left( x,y\right) $ of the CTSRW on $G$ satisfy 
\begin{equation*}
\sum_{x\in V(G)}\sum_{y\in \Gamma \left( x\right) }h(y,x)=n(n-1).
\end{equation*}
\end{lemma}

\begin{proof}
We know 
\begin{equation*}
h(x,x)=\frac{1}{d(x)}+\frac{1}{d(x)}\sum_{y\in \Gamma \left( x\right)
}h(y,x)=\frac{n}{d\left( x\right) }.
\end{equation*}
So, we have \ $\sum_{y\in \Gamma \left( x\right) }h(y,x)=n-1$ which is
independent of $x.$ Thus, for $n$ vertices, we will have 
\begin{equation*}
\sum_{x\in V(G)}\sum_{y\in \Gamma \left( x\right) }h(y,x)=n(n-1).
\end{equation*}
\end{proof}

If we denote the hitting time of $y$ from $x$ in SRW by $H(x,y)$,  then we
have an equality $$\sum_{x\in V(G)}\sum_{y\in \Gamma \left( x\right)
}H(y,x)=2m(n-1).$$
See \cite{Bo}. 
So, we have
the following theorem.

\begin{theorem}
With notations as the above, we have 
\begin{equation*}
\sum_{x\in V(G)}\sum_{y\in \Gamma \left( x\right) }H(y,x)=f\sum_{x\in
V(G)}\sum_{y\in \Gamma \left( x\right) }h(y,x).
\end{equation*}
\end{theorem}

Also, for CTSRW on a graph $G$, we define the mean commute time between vertices $x$
and $y$ to be $c(x,y)=h(x,y)+h(y,x)$. Let $C(x,y)$ be the corresponding quantity for SRW on $G$.
Then we have another version of the above equation in Theorem 2.4 as 
\begin{equation*}
\sum_{xy\in E(G)}C(x,y)=f\sum_{xy\in E(G)}c(x,y).
\end{equation*}

\section{Multi-person simple random walks on graphs}

We are led to the multi-person simple random walks (MPSRW) on a graph $G$ by the study of a continuous time Markov process induced by CTSRW on $G$. The combinatorics of MPSRW is much richer than we have touched upon here.    

Let $I_{n}$ be a finite set of cardinality $n$. For example, we may have   $I_{n}=\left\{1,2,\dots ,n\right\}$. We denote the set of all maps from $I_n$ to itself by $M_{n}$. We have the symmetric group $S_n$ sitting inside of $M_n$. A map $x\in M_{n}$ is called a {\it generalized permutation of deficiency} $k$ if 
\begin{equation*}
\left| x \left( I_{n}\right) \right| =n-k.
\end{equation*}
We denote $\text{def}(x)=k$.
\smallskip

$M_{n}$ is a semigroup under composition. The symmetrical group $S_{n}$
is a subgroup of $\ M_{n}.$ The deficiency determines a grading on \ $M_{n}$
which is compatible with the semigroup product on \ $M_{n}:$
\begin{equation*}
\text{def}\left(x\right) +\text{def}\left(y\right) \geq \text{def}\left( x\circ
y \right) .
\end{equation*}
Denote $M_{n}^{\left( k\right) }=\left\{ \text{all maps with deficiency }%
k\right\} ,$ then we have a decomposition of \ $M_{n}$ according the 
deficiency: 
\begin{equation*}
M_{n}=M_{n}^{\left( 0\right) }\amalg M_{n}^{\left( 1\right) }\amalg
M_{n}^{\left( 2\right) }\amalg \cdots \amalg M_{n}^{\left( n-1\right) },
\end{equation*}
where $M_{n}^{\left( 0\right) }=S_{n}.$ 
\smallskip

Let $G$ be a graph with the set of vertices $V(G)$ identified with $I_{n}$. We will call $G$ the {\it ground graph}.
Let the adjacency matrix $A\left( G\right)
=\left( a_{ij}\right) _{n\times n}$ with entries given by 
$$
a_{ij}=\begin{cases}
1 &\text{if $ij\in E\left( G\right)$}, \\ 
0 &\text{otherwise.}
\end{cases}
$$
\smallskip

We define a new graph $M\left( G\right)$ as follows: The set of vertices of 
$M\left( G\right)$ is $M_{n}$; For $x,y\in M_{n},$ there is an edge $xy$ in $M\left( G\right) $ only when 
\begin{equation*}
\left| \left\{ i:x\left( i\right) \neq y\left( i\right) \right\} \right| =1,
\end{equation*}
and if \ $x\left( i\right) \neq y\left( i\right) ,$ then \ $a_{x\left(
i\right) y\left( i\right) }=1.$ We will see that there is a close
relationship between the graph $M\left(G\right)$ and the $n$-th tensor power of CTSRW on the ground graph $G.$ 

\subsection{The tensor product of Markov processes}
Let $X_{t}^{\left( 1\right)},$ $X_{t}^{\left( 2\right) }$, $\dots
,X_{t}^{\left( n\right) }$ be Markov processes on the state spaces $
S^{\left( 1\right) },$ $S^{\left( 2\right) },\dots ,S^{\left( n\right) }$
respectively. We define a new process $Y_{t}$ on the state space 
$S^{\left( 1\right) }\times $ $S^{\left( 2\right)
}\times \cdots \times S^{\left( n\right) }$ with the transition probability 
given by 
$$\begin{aligned}
\Pr\,&\left\{ Y_{t+h}=\left( s_{2}\right) \mid Y_{h}=\left( s_{1}\right)
\right\} \\
&=\prod_{k=1}^{n}\Pr \left\{ Y_{t+h}^{\left( k\right) }=s_{2}^{\left(
k\right) }\mid Y_{h}^{\left( k\right) }=s_{1}^{\left( k\right) }\right\} \\
&=\prod _{k=1}^{n}p_{s_{1}^{\left( k\right) }s_{2}^{\left( k\right)
}}^{\left( k\right) }\left( t\right) ,
\end{aligned}
$$
where $\left( s_{i}\right) =\left( s_{i}^{\left( 1\right) },s_{i}^{\left( 2\right) },\dots ,s_{i}^{\left( n\right)
}\right)$, $i=1,2,$ and $p_{s_{1}^{\left( k\right) }s_{2}^{\left(
k\right) }}^{\left( k\right) }\left( t\right) $ is the transition
probability of the Markov process $X_{t}^{\left( k\right) }$, $k=1,2,\dots,n.$
\smallskip

We call $Y_{t}$ the {\it tensor product} of Markov processes $X_{t}^{\left(
k\right) },$ $k=1,2,\cdots ,n.$
\smallskip

The next two lemmas can be proved by some direct computations. So we omit the proofs.
 
\begin{lemma}
$Y_{t}$ is a Markov process.
\end{lemma}

\begin{lemma}
Let $q_{\left( s_{1}\right) \left( s_{2}\right) }$ be the infinitesmal generator
of $Y_{t}$. Then 
$$
q_{\left( s_{1}\right) \left( s_{2}\right) }=\begin{cases}
q_{s_{1}^{\left( k\right) }s_{2}^{\left( k\right) }}^{\left( k\right)}
&\text{if $\exists$ only one index $k$ such that $s_{1}^{\left( k\right)}\neq s_{2}^{\left( k\right)}$},\\ 
\sum_{k=1}^{n}q_{s_{1}^{\left( k\right) }s_{2}^{\left( k\right) }}^{\left(
k\right) }&\text{if $\left(s_{1}\right) =\left( s_{2}\right)$, namely, $s_{1}^{\left( k\right)}=s_{2}^{\left( k\right)}$ for all $k$}, \\ 
0&\text{otherwise.}
\end{cases}
$$
where $q_{s_{1}^{\left( k\right) }s_{2}^{\left( k\right) }}^{\left( k\right)
}$ is the infinitesmal generator of Markov process \ $X_{t}^{\left( k\right)
},$ $k=1,2,\dots ,n.$
\end{lemma}

Now, if we order the elements of the state space $S^{\left( 1\right) }\times
S^{\left( 2\right) }\times \cdots \times S^{\left( n\right) }$
lexicographically, it is easy to see that the probability transition matrix
of $Y_{t}$ is given by the tensor product 
\begin{equation*}
P^{\left( 1\right) }\left( t\right) \otimes P^{\left( 2\right) }\left(
t\right) \otimes \cdots \otimes P^{\left( n\right) }\left( t\right).
\end{equation*}
This is why we call $Y_{t}$ is the tensor product of Markov process $%
X^{\left( k\right) }\left( t\right),$ $k=1,2,\dots ,n.$ By one of the lemmas above, we also can see that the infinitesimal generator matrix of $Y_{t}$
is given by 
\begin{equation*}
Q^{\left( 1\right) }\otimes I\otimes \cdots \otimes I+I\otimes Q^{\left(
2\right) }\otimes I\otimes \cdots \otimes I+\cdots \cdots +I\otimes \cdots
\otimes I\otimes Q^{\left( n\right) }
\end{equation*}
where $Q^{\left( k\right) }$ is the infinitesimal generator matrix of $%
X_{t}^{\left( k\right) },k=1,2,\dots n.$ Thus, for convenience, we denote 
$Y_{t}=X_{t}^{\left( 1\right) }\otimes X_{t}^{\left( 2\right) }\otimes
\cdots \otimes X_{t}^{\left( n\right) }$.

\subsection{CTSRW on M(G)}
On the graph $M(G)$, we have the Markov process CTSRW. Denote it by $Y_{t}$. 
Let $X_t$ be the Markov process CTSRW on $G$. 

\begin{lemma}
The Markov process $Y_{t}$ on the graph $M(G)$ is the $n$-th 
tensor power of the Markov process $X_{t}$ on the ground graph $G.$
The jump chain of $Y_t$ is MPSRW of $G$.
\end{lemma}

\begin{proof}
Let's recall $Q=-L\left( G\right) =-D\left( G\right) +A\left( G\right)
=\left( q_{ij}\right) _{n\times n},$%
\begin{equation*}
q_{ij}=\begin{cases}
1&\text{if $ij\in E\left( G\right)$,} \\ 
-d\left( i\right)&\text{if $i=j$,} \\ 
0&\text{otherwise,}
\end{cases}
\end{equation*}
where $d\left( i\right) $ is the degree of the vertex $i.$ Let $X_{t}\otimes
X_{t}\otimes \cdots \otimes X_{t}=Z_{t}$ be the $n$-th tensor power of
the process $X_{t}$. Take two states for $Z_{t}$, 
$\left( s_{1}\right) $ and $%
\left( s_{2}\right)$. Then $\left( s_{i}\right) $ actually is a sequence
of vertices of $G.$ We write $\left( s_{1}\right) =\left(
i_{1},i_{2},\dots ,i_{n}\right) $ and $\left( s_{2}\right) =\left(
j_{1},j_{2},\dots ,j_{n}\right).$ By the lemma above, we have
$$\begin{aligned}
p_{\left( s_{1}\right) \left( s_{2}\right) }^{\prime }\left( 0\right) 
&=\begin{cases}
q_{s_{1}^{\left( k\right) }s_{2}^{\left( k\right) }}^{\left( k\right) } &
\text{if $\exists$ only one index $k$ such that $s_{1}^{\left( k\right) }\neq
s_{2}^{\left( k\right) }$,} \\ 
\sum_{k=1}^{n}q_{s_{1}^{\left( k\right) }s_{2}^{\left( k\right) }}^{\left(
k\right) }&\text{if $\left( s_{1}\right) =\left( s_{2}\right)$, 
namely, $s_{1}^{\left( k\right) }=s_{2}^{\left( k\right) }$ for all $k$}, \\ 
0&\text{otherwise}
\end{cases}\\
&=\begin{cases}
1&\text{if $\exists$ only one index $k$ such that  $i_{k}\neq j_{k},i_{k}j_{k}\in E\left(
G\right)$} , \\ 
-\sum_{k=1}^{n}d\left( i_{k}\right)&\text{if $i_{k}=j_{k}$ for all $k$},\\ 
0&\text{otherwise.}
\end{cases}
\end{aligned}
$$

Now, we can identify $\left( s_{1}\right)$ and $(s_2)$ with the images of certain maps $x$ and $y$, respectively.
Then by the definition of $M(G)$, $xy\in E\left( M\left( G\right) \right) ,
$ if and only if $p_{(s_1)(s_2)}^{\prime }(0)=1.$ As to
the degree of $x$, we know that the neighbors of $x=\left( s_{1}\right) $
can only be $\left( k_{1},i_{2},\dots ,i_{n}\right)$, where $k_{1}$ must
be a neighbor of $i_{1}$; $\left( i_{1},k_{2},\dots ,i_{n}\right) ,$ where 
$k_{2}$ must be a neighbor of $i_{2};$ and so on, up to the last one $%
\left( i_{1},i_{2},\dots ,i_{n-1},k_{n}\right) ,$ where $k_{n}$ must be a
neighbor of $i_{n}$. Thus 
$$d(x)=d(s_{1})=\sum_{j=1}^{n}d\left(
i_{j}\right) .$$ 
This agree with $p_{\left( s_{1}\right) \left( s_{1}\right)
}^{\prime }(0)=-\sum_{k=1}^{n}d\left( i_{k}\right) .$ \ So, \ $Y_{t}=
$\ $X_{t}\otimes X_{t}\otimes \cdots \otimes X_{t}=Z_{t}.$
\smallskip

The second conclusion is easy to see.
\end{proof}

Because of this lemma, we may call graph $M(G)$ the tensor power of the graph $
G.$ We also note that the transition probability from $x$ to $y$ in $M_n$ is given by  
\begin{equation*}
P_{xy}\left( t\right) =\prod_{i=1}^{n}p_{x\left( i\right) y\left( i\right)
}\left( t\right) .
\end{equation*}
\smallskip

From the proof of Lemma 3.3, we know that for $x\in M_{n}=V(M(G))$, 
$$d(x)=\sum_{k=1}^{n}d(x(k)).$$ We can compute the
size of $M\left( G\right)$, i. e. the number of edges of $M(G)$ as 
\begin{eqnarray*}
2\left| E(M(G))\right| &=&\sum_{x\in M_{n}}d\left( x\right) =\sum_{x\in
M_{n}}\sum_{k=1}^{n}d\left( x\left( k\right) \right)
=\sum_{k=1}^{n}\sum_{x\in M_{n}}d\left( x\left( k\right) \right) \\
&=&\sum_{k=1}^{n}\left( \left( d\left( 1\right) +d\left( 2\right) +\cdots
+d\left( n\right) \right) n^{n-1}\right) \\
&=&2mn^{n}.
\end{eqnarray*}
Thus, the size of $M(G)$ is $mn^{n}.$ It is clear that the order of $M(G)$ is $n^n$. Therefore, the frequency of $M(G)$ is $2m$.

\begin{lemma}
$G$ is bipartite if and only if $M\left( G\right) $ is bipartite.
\end{lemma}

\begin{proof}
We use the classical result of K\"{o}nig that a graph is bipartite if and only if all its cycles are even.
\smallskip

If $G$ is not bipartite, then there is a cycle $C=
v_{1}v_{2}\cdots v_{m}v_{1}$ of odd length in $G$. We look
at the cycle in $M\left( G\right)$ given by
\begin{equation*}
\left( v_{1},v_{1},\cdots ,v_{1}\right)\left( v_{1},\cdots
,v_{1},v_{2}\right)\left( v_{1},\cdots ,v_{1},v_{3}\right)\cdots
\left( v_{1},\cdots ,v_{1},v_{m}\right)\left( v_{1},\cdots
,v_{1},v_{1}\right). 
\end{equation*}
Its length is also odd. Thus $M(G)$ is not bipartite.
\smallskip

If $G$ is bipartite, then the set of vertices $V=V(G)$ can be written as 
$V_{1}\cup V_{2}$, with $V_{1}\cap V_{2}=\emptyset$ and there is no edge between vertices both in $V_1$ or both in $V_2$. We try to bipart the
set of vertices of $M(G).$ We know $V(M(G))=V^{\times n}$. We write $V=V_{1}+V_{2}$. Then 
\begin{equation*}
V^{\times n}=\left( V_{1}\cup V_{2}\right) ^{\times n}=
=V_{1}^{n}\cup V_{1}^{n-1}V_{2}\cup V_{1}^{n-2}V_{2}^{2}\cup\dots\cup
V_{1}V_{2}^{n-1}\cup V_{2}^{n},
\end{equation*}
where $V_{1}^{n-k}V_{2}^{k}$ means we take $\left( n-k\right) $
vertices from $V_{1}$ and $k$ vertices from $V_{2}$, regardless of order, to form a vertex of $M(G).$ Let 
\begin{eqnarray*}
\overline{V_{1}}&=&V_{1}^{n}\cup V_{1}^{n-2}V_{2}^{2}\cup
V_{1}^{n-4}V_{2}^{4}\cup\cdots, \\
\overline{V_{2}}&=&V_{1}^{n-1}V_{2}\cup V_{1}^{n-3}V_{2}^{3}\cup V_{1}^{5}V_{2}^{n-5}\cup\cdots.
\end{eqnarray*}
Then $\overline{V_{1}}\cup\overline{V_{2}}=V(M(G))$ and $\overline{
V_{1}}\cap \overline{V_{2}}=\emptyset .$ By the definition of $M(G)$, we
can not find an edge between any two vertices which are both in $\overline{%
V_{1}}$ or both in $\overline{V_{2}}.$
\end{proof}

If $G$ is not bipartite with order $n$ and size $m.$ Write the degree sequence of $G$ as $d_{1}\leq d_{2}\leq \cdots \leq d_{n}.$ For
any vertex $x$ in $M(G),$ the $x$-component of the ergodic vector for SRW
on $M(G)$ is given by 
\begin{equation*}
\pi_{x}=\frac{d(x)}{2mn^{n}}=\frac{\sum_{i=1}^{n}d(x\left(i\right))}{
2mn^{n}}=n^{-n}\sum_{i=1}^{n}\frac{d(x\left( i\right) )}{2m}%
=n^{-n}\sum_{i=1}^{n}\pi _{x\left( i\right) }^{G},
\end{equation*}
where 
$$\pi _{x\left( i\right) }^{G}=\frac{d(x\left( i\right))}{2m}$$ 
is the $x(i)$-component of the ergodic vector for SRW on $G.$
\smallskip

The expected first return time for CTSRW on $M(G)$ is given by 
\begin{equation*}
h(x,x)=\frac{n^{n}}{d\left( x\right) }=\frac{n^{n}}{\sum_{k=1}^{n}d(x%
\left( k\right))}.
\end{equation*}
In particular, for the identity map $Id$, $d(Id)=d_{1}+d_{2}+\cdots +d_{n}=2m$. So 
\begin{equation*}
h\left(Id,Id\right) =\frac{n^{n}}{2m}.
\end{equation*}
\smallskip

For the jump chain SRW on $M(G)$, or MPSRW on $G$, we have
\begin{equation*}
H(x,x)=f_{M(G)}h\left( x,x\right) =\frac{2mn^{n}}{\sum_{k=1}^{n}
d(x(k))},
\end{equation*}
where $f_{M(G)}=2m$ is the frequency of $M(G)$. In particular,
\begin{equation*}
H\left( Id,Id\right) =n^{n}.
\end{equation*}

\subsection{Hitting time for CTSRW on $M(G)$} For CTSRW on $M(G)$, we are interested in calculating $h(Id,c_i)$, where $c_i$ is the constant map to a vertex $i$ of $G$. Once we know $h(Id,i)$, we can use it to estimate $H(Id,c_i)$ for SRW on $M(G)$, or equivalently, the expect number of steps MPSRW on $G$ should take to have all persons meet at the vertex $i$. Namely, we have
\begin{equation*}
nd_{1}h(Id,c_{i}) \leq H(Id,c_{i}) \leq
nd_{n}h(Id,c_{i}).
\end{equation*}
Notice that when the graph $G$ is regular, $H(Id,c_i)$ is determined completely by $h(Id,c_i)$: $H(Id,c_i)=nd\cdot h(Id,c_i)$. 
\smallskip

Let $x$ and $y$ be two distinct vertices in $M(G)$, we consider the hitting
time $h(x,y)$ of $y$ from $x$ in CTSRW. We have 
\begin{eqnarray*}
h\left( x,y\right) &=&n^{n}\int_{0}^{\infty }\left( P_{yy}\left( t\right)
-P_{xy}\left( t\right) \right) dt \\
&=&n^{n}\int_{0}^{\infty }\left( \prod_{i=1}^{n}p_{y\left( i\right) y\left(
i\right) }\left( t\right) -\prod_{i=1}^{n}p_{x\left( i\right) y\left(
i\right) }\left( t\right) \right) dt.
\end{eqnarray*}
For $x=Id,$ $y=c_{i}$, we have 
\begin{equation*}
h\left(Id,c_{i}\right) =n^{n}\int_{0}^{\infty }\left( p_{ii}(t)^{n}-\prod_{j=1}^{n}p_{ji}(t) \right) dt.
\end{equation*}
\smallskip

To calculate this integral, we diagonalize the Laplacian of $G$. Write 
$$U^{T}QU=\text{diag}\left[ -\lambda _{1},-\lambda _{2},\cdots ,-\lambda _{n}\right],$$ 
where $U=\left(u_{ij}\right) _{n\times n}$ is an orthogonal matrix. Then 
\begin{equation*}
p_{ij}\left( t\right) =\sum_{k=1}^{n}u_{ik}u_{jk}e^{-\lambda _{k}t}.
\end{equation*}
Since 
\begin{equation*}
p_{ii}(t)^{n}=\left( \sum_{k=1}^{n}u_{ik}^{2}e^{-\lambda
_{k}t}\right)^{n}=\sum_{1\leq k_{1},k_{2},\cdots ,k_{n}\leq
n}u_{ik_{1}}^{2}u_{ik_{2}}^{2}\cdots u_{ik_{n}}^{2}e^{-(\lambda
_{k_{1}}+\lambda _{k_{2}}+\cdots +\lambda _{k_{n}})\, t}
\end{equation*}
and
\begin{equation*}
\prod_{k=1}^{n}p_{ki}\left( t\right) =\sum_{1\leq k_{1},k_{2},\cdots
,k_{n}\leq n}u_{1k_{1}}u_{ik_{1}}u_{2k_{2}}u_{ik_{2}}\cdots
u_{nk_{n}}u_{ik_{n}}e^{-\left( \lambda _{k_{1}}+\lambda _{k_{2}}+\cdots
+\lambda _{k_{n}}\right)\,t},
\end{equation*}
we get 
\begin{equation*}
h\left(Id,c_{i}\right) =n^n
\sum_{1\leq k_{1},k_{2},\cdots ,k_{n}\leq n}
\frac{u_{ik_{1}}^{2}u_{ik_{2}}^{2}\cdots u_{ik_{n}}^{2}-
u_{1k_{1}}u_{ik_{1}}u_{2k_{2}}u_{ik_{2}}\cdots u_{nk_{n}}u_{ik_{n}}}
{\lambda_{k_{1}}+\lambda _{k_{2}}+\cdots +\lambda _{k_{n}} }.
\end{equation*}
\smallskip

\noindent{\bf Example.} Let $G$ be the triangle graph. It is regular 
with $n=3$ and $d=2$. The matrix $Q$ and $U$ are given below:
$$
Q=\begin{pmatrix}
-2 & 1 & 1 \\
1 & -2 &1 \\
1 & 1 & -2 
\end{pmatrix}
\quad\text{and}\quad
U=\begin{pmatrix}
-\frac1{\sqrt{2}} & -\frac1{\sqrt{6}} & \frac1{\sqrt{3}} \\
0 & \frac 2{\sqrt{6}} & \frac1{\sqrt{3}} \\
\frac1{\sqrt{2}} & -\frac1{\sqrt{6}} & \frac1{\sqrt{3}}
\end{pmatrix}.
$$
We have $U^TQU=\text{diag}[-3,-3,0]$. Since $\lambda_3=0$, when 
calculating $h(Id,c_i)$, we should drop the term $(k_1,k_2,k_3)=(3,3,3)$. 
Notice that the numerator of this term in $h(Id,c_i)$ is also zero, so it is 
fine to drop this term. We have
$$h(Id,c_1)=3^3\cdot\frac{31}{162}=\frac{31}{6}.$$
So 
$$H(Id,c_1)=3\cdot2\cdot\frac{31}{6}=31.$$
Thus, for our MPSRW on the triangle graph, it takes 31 steps on average 
for 3 persons to meet at any specified vertex, given that they all 
start at different 
vertices.
\smallskip

For the square graph with $n=4$ and $d=2$, a similar calculation shows
$$h(Id,c_1)=4^4\cdot\frac {167}{1120}=\frac{1336}{35}$$
and 
$$H(Id,c_1)=8\cdot\frac{1336}{35}\approx 305.371.$$
Thus, for our MPSRW on the square graph, it takes about 305 steps on 
average for 4 persons to meet at any specified vertex, given that they all start
at different vertices.

\end{document}